\documentclass[aps,onecolumn]{revtex4}
\usepackage[dvips]{graphicx}
\newcommand{\sect}[1]{\setcounter{equation}{0}\section{#1}}

\begin{document}
\title{Theoretical foundations and mathematical formalism of the power-law tailed  statistical distributions}
\author{G. Kaniadakis}
\email{giorgio.kaniadakis@polito.it}
\affiliation{Department of Applied Science and Tecnology, Politecnico di Torino, \\
Corso Duca degli Abruzzi 24, 10129 Torino, Italy}
\date{\today}
\begin {abstract}
We present the main features of  the mathematical theory generated by the $\kappa$-deformed exponential function $\exp_{\kappa}(x)=(\sqrt{1+ \kappa^2 x^2}+\kappa x)^{1/\kappa}$, with $0\leq \kappa<1$, developed in the last twelve years, which turns out to be a continuous one parameter deformation of the ordinary mathematics generated by the Euler exponential function.  The $\kappa$-mathematics has its roots in special relativity and furnishes the theoretical foundations of the $\kappa$-statistical mechanics predicting power law tailed statistical distributions which have been observed experimentally in many physical, natural and artificial systems. After introducing the $\kappa$-algebra
we present the associated $\kappa$-differential and $\kappa$-integral calculus. Then we obtain the corresponding $\kappa$-exponential and $\kappa$-logarithm functions and give the $\kappa$-version of the main functions of the ordinary mathematics.

\end {abstract}


\maketitle

\sect{Introduction}

Undoubtedly the most interesting feature of the statistical distribution function
\begin{equation}
f_i=\exp_{\kappa}(-\beta E_i+\beta \mu) \ \ ,  \label{A1}
\end{equation}
where the $\kappa$-exponential is defined as
\begin{equation}
\exp_{\kappa}(x)=(\sqrt{1+ \kappa^2 x^2}+\kappa x)^{1/\kappa} \ \ , \ \ \ 0\leq\kappa<1 \ \ ,
\label{A2}
\end{equation}
is represented by its asymptotic behavior, namely
\begin{equation}
f_i{\atop\stackrel{\textstyle\sim}{\scriptstyle {\beta E_i-\beta\mu} \rightarrow
\,0}} \exp(-\beta E_i+\beta \mu) \ \ , \ \ \ \ f_i {\atop\stackrel{\textstyle\sim}{\scriptstyle E_i \rightarrow
\,+\infty}} N E_i^{-1/\kappa} \ \ .   \label{A3}
\end{equation}
The above $\kappa$-distribution at low energies is the ordinary Boltzmann distribution while at high energies presents an power-law  tail. For this reason the statistical theory \cite{PhA2001,PLA2001,PRE2002,PRE2005,PhA2006,EPJA2009,EPJB2009,EPL2010,PLA2011,MPLB2012},
based on the distribution (\ref{A1}) has attracted the interest of many researchers.

In the last twelve years various authors have considered the
foundations of the statistical theory based on the $\kappa$-distribution, in connection with the historical evolution of the research on the power-law tailed statistical distributions \cite{EPJB2009,KL2004}
e.g. the H-theorem and the molecular chaos hypothesis
\cite{Silva06A,Silva06B}, the thermodynamic stability
\cite{Wada1,Wada2}, the Lesche stability
\cite{KSPA04,AKSJPA04,Naudts1,Naudts2}, the Legendre structure of
the ensued thermodynamics \cite{ScarfoneWada,Yamano},
the thermodynamics of non-equilibrium systems \cite{Lucia2010},
quantum versions of the theory \cite{AlianoKM2003,Santos2011a,Santos2011b,Santos2012},
the geometrical structure of the theory \cite{Pistone},
various mathematical aspects of the theory
\cite{KLS2004,KLS2005,KSphysA2003,Oikonomou2010,Stankovic2011,
Tempesta2011,DeossaCasas,Vigelis,Scarfone2013}, etc.
On the other hand specific applications to physical systems have been
considered, e.g. the cosmic rays \cite{PRE2002}, relativistic
\cite{GuoRelativistic} and classical \cite{GuoClassic} plasmas  in
presence of external electromagnetic fields, the relaxation in
relativistic plasmas under wave-particle interactions
\cite{Lapenta,Lapenta2009},  anomalous diffusion
\cite{WadaScarfone2009,Wada2010}, non-linear kinetics \cite{KQSphysA2003,BiroK2006,Casas},
kinetics of interacting atoms and
photons \cite{Rossani}, particle kinetics in the presence of
temperature gradients \cite{GuoDuoTgradient,Guo2012},  particle systems in
external conservative force fields \cite{Silva2008}, stellar
distributions in astrophysics \cite{Carvalho,Carvalho2,Carvalho2010,Bento2013},
quark-gluon plasma formation \cite{Tewel}, quantum hadrodynamics models
\cite{Pereira}, the fracture propagation \cite{Fracture}, etc. Other
applications concern dynamical systems at the edge of chaos
\cite{Corradu,Tonelli,Celikoglu}, fractal systems \cite{Olemskoi},
field theories \cite{Olemskoi2010},
the random matrix theory \cite{AbulMagd,AbulMagd2009,AbulMagd2012}, the error
theory \cite{WadaSuyari06}, the game theory \cite{Topsoe}, the theory of complex networks \cite{Macedo}, the
information theory \cite{WadaSuyari07}, etc. Also applications to
economic systems have been considered e.g. to study the personal
income distribution \cite{Clementi,Clementi2008,Clementi2009,Clementi2011,Clementi2012a,Clementi2012b}, to model deterministic heterogeneity in tastes and product
differentiation \cite{Rajaon,Rajaon2008}, in finance \cite{Trivellato2012,Trivellato2013}, in equity options \cite{Tapiero}, to construct  taxation and redistribution models \cite {Bertotti}, etc.

In this contribution we present the main features of  the mathematical theory generated by the function $\exp_{\kappa}(x)$. The  $\kappa$-mathematics, developed in the last twelve years, turns out to be a continuous one parameter deformation of the ordinary mathematics generated by the Euler exponential function.  The $\kappa$-mathematics has its roots in special relativity and furnishes the theoretical foundations of the $\kappa$-statistical theory predicting power law tailed statistical distributions which have been observed experimentally in many physical, natural and artificial systems.

The paper is organized as follows: After introducing the $\kappa$-algebra
we present the associated $\kappa$-differential and $\kappa$-integral calculus. Then we obtain the corresponding $\kappa$-generalized exponential and logarithmic functions and give the $\kappa$-version of the main functions of the ordinary mathematics.

\sect{$\kappa$-Algebra}

\noindent {\bf Theorem II.1.} Let be $x,y \in {\bf R}$ and $-1<\kappa<1$. The composition law $\stackrel{\kappa}{\oplus}$ defined through
\begin{equation}
x\stackrel{\kappa}{\oplus}y=x\sqrt{1+\kappa^2y^2}+y\sqrt{1+\kappa^2x^2}
\ \ \nonumber , \label{MII1}
\end{equation}
is a generalized sum, called $\kappa$-sum and the algebraic structure $({\bf
R},\stackrel{\kappa}{\oplus})$ forms an abelian group.

\noindent {\bf Proof.} From the definition of $\stackrel{\kappa}{\oplus}$ the
following properties follow
\\ \noindent 1) associativity:
$(x\stackrel{\kappa}{\oplus}y)\stackrel{\kappa}{\oplus}z
=x\stackrel{\kappa}{\oplus}(y\stackrel{\kappa}{\oplus}z)$,
\\ \noindent 2) neutral element: $x\stackrel{\kappa}{\oplus}0=0
\stackrel{\kappa}{\oplus}x=x$, \\ \noindent 3) opposite element:
$x\stackrel{\kappa}{\oplus}(-x)=(-x) \stackrel{\kappa}{\oplus}x=0$,
\\ \noindent 4) commutativity:
$x\stackrel{\kappa}{\oplus}y=y\stackrel{\kappa}{\oplus}x$.

\noindent {\bf Remarks.} The $\kappa$-sum is a one parameter continuous deformation of the ordinary sum which recovers in the classical limit
$\kappa\rightarrow 0$, i.e. $x\stackrel{0}{\oplus}y=x+y$.  The $\kappa$-sum (\ref{MII1}) is the additivity law of the dimensionless relativistic momenta of special relativity while the real parameter  $-1<\kappa<1$ is the reciprocal of the dimensionless light speed \cite{PRE2002,PLA2011}. The $\kappa$-difference
$\stackrel{\kappa}{\ominus}$ is defined as
$x\stackrel{\kappa}{\ominus}y=x\stackrel{\kappa}{\oplus}(-y)$.

\noindent {\bf Theorem II.2.} Let be $x,y \in {\bf R}$ and $-1<\kappa<1$. The composition law $\stackrel{\kappa}{\otimes}$
defined through
\begin{equation}
x\stackrel{\kappa}{\otimes}y={1\over\kappa}\,\sinh
\left(\,{1\over\kappa}\,\,{\rm arcsinh}\,(\kappa x)\,\,{\rm
arcsinh}\,(\kappa y)\,\right) \ , \label{MII2}
\end{equation}
is a generalized product, called $\kappa$-product and
the algebraic structure $({\bf R},\stackrel{\kappa}{\otimes})$
forms an abelian group.

\noindent {\bf Proof.} From the definition of $\stackrel{\kappa}{\otimes}$ the
following properties follow
\\ \noindent 1) associativity:  $(x \stackrel{\kappa}{\otimes}y)
 \stackrel{\kappa}{\otimes}z=x
 \stackrel{\kappa}{\otimes}(y
 \stackrel{\kappa}{\otimes}z)$,
\\ \noindent 2) neutral element: is defined through $x \stackrel{\kappa}{\otimes}I=I
 \stackrel{\kappa}{\otimes}x= x$ and is given by $I=\kappa^{-1}\sinh
\kappa$,
\\ \noindent 3) inverse element: is defined through  $x \stackrel{\kappa}{\otimes}\overline x= \overline x
\stackrel{\kappa}{\otimes}x=I$ and is given by $\overline
x=\kappa^{-1}\sinh(\kappa^2/{\rm arcsinh} \,\kappa x)$,
\\ \noindent 4) commutativity: $x\stackrel{\kappa}{\otimes}y= y
 \stackrel{\kappa}{\otimes}x$.

\noindent {\bf Remarks.} The $\kappa$-product reduces to the ordinary product as $\kappa\rightarrow 0$, i.e. $x\stackrel{0}{\otimes}y=xy$. The $\kappa$-division
$\stackrel{\kappa}{\oslash}$ is defined through
$x\stackrel{\kappa}{\oslash}y=x\stackrel{\kappa}{\otimes}\overline
y$.

\noindent {\bf Theorem II.3.} Let be $x,y \in {\bf R}$ and $-1<\kappa<1$. The $\kappa$-sum $\stackrel{\kappa}{\oplus}$ defined in (\ref{MII1}), and the $\kappa$-product $\stackrel{\kappa}{\otimes}$ defined in (\ref{MII2}), obey the distributive law
\begin{equation}
z\stackrel{\kappa}{\otimes}(x \stackrel{\kappa}{\oplus}y) = (z
\stackrel{\kappa}{\otimes}x) \stackrel{\kappa}{\oplus}(z
\stackrel{\kappa}{\otimes}y) \ \ , \label{MII3}
\end{equation}
and then the algebraic structure $({\bf R},\stackrel{\kappa}
{\oplus},\stackrel{\kappa}{\otimes})$ forms an abelian field.

\noindent {\bf Proof.} The relationship (\ref{MII3}) follows directly from the definitions
of the $\kappa$-product (\ref{MII2}) and of the $\kappa$-sum
(\ref{MII1}) which can be written also in the form
\begin{equation}
x\stackrel{\kappa}{\oplus}y={1\over\kappa}\,\sinh\Big(\,{\rm
arcsinh}\,(\kappa x)+{\rm arcsinh}\,(\kappa y)\,\Big) \ .
\label{MII4}
\end{equation}

\noindent {\bf Theorem II.4.} The abelian fields $({\bf R},\stackrel{\kappa}
{\oplus},\stackrel{\kappa}{\otimes})$ and $({\bf R}, + , .)$ are
isomorphic.

\noindent {\bf Proof.} After introducing the function $\{x\}\in
C^{\infty}({\bf R})$ through
\begin{equation}
\{x\}=\frac{1}{\kappa}\,{\rm arcsinh} \,(\kappa x) \ \ ,
\label{MII5}
\end{equation}
whose inverse function $[x]\in C^{\infty}({\bf R})$, i.e.
$[\{x\}\,]=\{\,[x]\}=x$, is given by
\begin{equation}
[x]=\frac{1}{\kappa}\,\sinh \,(\kappa x) \ \ , \label{MII6}
\end{equation}
we can write Eqs. (\ref{MII4}) and (\ref{MII2}) in the form
\begin{eqnarray}
&&\{x\stackrel{\kappa}{\oplus}y\}=\{x\}+\{y\} \ , \label{MII7}
\\ &&\{x\stackrel{\kappa}{\otimes}y\}=\{x\}\cdot\{y\}  \ , \label{MII8}
\end{eqnarray}
or equivalently as
\begin{eqnarray}
&&[x]\stackrel{\kappa}{\oplus}[y]=[x+y] \ , \label{MII9}
\\ &&[x]\stackrel{\kappa}{\otimes}[y]=[x\cdot y]  \ . \label{MII10}
\end{eqnarray}

\noindent {\bf Theorem II.5.} Let be $x \in {\bf R}$ and $n$ an arbitrary nonnegative integer.  It holds
\begin{equation}
 \underbrace{x\stackrel{\kappa}{\oplus}x\stackrel{\kappa}{\oplus}...\stackrel{\kappa}{\oplus}x
}_{n \ times}= [n]\stackrel{\kappa}{\otimes}
 x  \  . \ \
\label{MII11}
\end{equation}
\noindent {\bf Proof.} The function $[x]$ and its inverse $\{x\}$ obey the condition $[\{x\}]=\{[x]\}=x$. Furthermore we take into account (\ref{MII7}) and (\ref{MII10}). Then we have
\begin{eqnarray}
x\stackrel{\kappa}{\oplus}x\stackrel{\kappa}{\oplus}...\stackrel{\kappa}{\oplus}x
&=& [\{x\stackrel{\kappa}{\oplus}x\stackrel{\kappa}{\oplus}...\stackrel{\kappa}{\oplus}x \}] \nonumber
\\ &=& [\{x\}+\{x\}+ ... + \{x\}] \nonumber
\\ \nonumber &=& [n\cdot\{x\}] \\ \nonumber &=& [\{[n]\}\cdot\{x\}] \\ \nonumber &=&[n]\stackrel{\kappa}{\otimes}
[\{x\}]\\ \nonumber &=&[n]\stackrel{\kappa}{\otimes}x \nonumber \ \ .
\end{eqnarray}

\sect{$\kappa$-Differential Calculus}

\subsection{$\kappa$-Differential}

The $\kappa$-differential of $x$, indicated
by  $d _{\kappa}x$, is defined through
\begin{eqnarray}
(x+dx) \stackrel{\kappa}{\ominus}x=d _{\kappa}x + 0((dx)^2) \ , \label{MIII1}
\end{eqnarray}
and results to be
\begin{eqnarray}
 d _{\kappa}x= \frac{d\,x}{\displaystyle{\sqrt{1+\kappa^2\,x^2} }} \ .
\label{MIII2}
\end{eqnarray}

In order to better understand the origin of the expression of $d
_{\kappa}x$ we recall that the variable $x$ is a dimensionless
momentum. Then the quantity $\gamma(x)=\sqrt{1+\kappa^2\,x^2}$ is
the Lorentz factor of relativistic physics, in the momentum representation. So, we can write
\begin{eqnarray}
 d _{\kappa}x= \frac{dx}{\gamma(x)} \ .
\label{MIII2a}
\end{eqnarray}

Moreover it holds
\begin{eqnarray}
 d _{\kappa}x=d \{x\}= \frac{d \{x\}}{dx}\,\,dx \ . \label{MIII3}
\end{eqnarray}

\subsection{$\kappa$-Derivative}

We define the $\kappa$-derivative of the
function $f(x)$ through
\begin{equation}
\frac{d\,f(x)}{d_{\kappa} \,x}=\lim_{z\rightarrow x}\frac{f
(z)-f(x)}{\displaystyle{ z \stackrel{\kappa}{\ominus}x}}\approx \lim_{dx\rightarrow 0}\frac{f
(x+dx)-f(x)}{\displaystyle{ (x+dx) \stackrel{\kappa}{\ominus}x}}\
\ . \label{MIII4}
\end{equation}

We observe that $df(x)/d_{\kappa}x$, which reduces to $df(x)/dx$
as the deformation parameter ${\kappa}\rightarrow 0$, can be
written in the form
\begin{eqnarray}
\frac{d \, f(x)}{d _{\kappa}\,x}=\sqrt{1+\kappa^2\,x^2}\,\,\,\frac{d
\, f(x)}{d \, x} \ . \label{MIII5}
\end{eqnarray}
From the latter relatioship follows that the $\kappa$-derivative
obeys the Leibniz's rules of the ordinary derivative. After
introducing the $\gamma(x)$ Lorentz factor the $\kappa$-derivative
can be written also in the form:
\begin{eqnarray}
\frac{d }{d _{\kappa}\,x}=\gamma(x)\,\frac{d }{d \, x} \ .
\label{MIII5a}
\end{eqnarray}

\subsection{$\kappa$-Integral}

We define the $\kappa$-integral as the inverse operator of the $\kappa$-derivative through
\begin{eqnarray}
\int d_{\kappa}x \,\, f(x)= \int \frac{d \,
x}{\sqrt{1+\kappa^2\,x^2}}\,\,f(x) \ , \label{MIII6}
\end{eqnarray}
and note that it is governed by the same rules of the ordinary
integral, which recovers when  ${\kappa}\rightarrow 0$.

\subsection{Connections with Physics}

We indicate with $p=|{\bf p}|$ and $x=p/mv_*$ the moduli of the
particle momentum in dimensional and dimensionless form respectively and define
$\kappa=v_*/c$. The classical relationship linking $x$ with the
dimensionless kinetic energy ${\cal W}=x^2/2$ follows from the
kinetic energy theorem, which in differential form reads
\begin{eqnarray}
\frac{d}{dx}\, {\cal W} =x \ . \label{MIII7}
\end{eqnarray}
The latter equation after replacing the ordinary derivative by the derivative $d/d_{\kappa}x$, i.e.
\begin{eqnarray}
\frac{d}{d_{\kappa}x}\, {\cal W}=x  \ , \label{MIII8}
\end{eqnarray}
transforms into the corresponding relativistic equation. This differential equation with the condition
${\cal W}(x=0)=0$ admits as unique solution
${\cal W}=\left(\sqrt{1+\kappa^2x^2}-1\right)/\kappa^2$ defining  the relativistic  kinetic energy.

Let us consider the four-dimensional Lorentz invariant integral
\begin{eqnarray}
I = \int d^4\!p \,\,\theta (p_0)\,
\delta(p^{\mu}p_{\mu}-m^2c^2)\,F(p) \ , \label{MIII9}
\end{eqnarray}
being $p^{\mu}=(p^0,{\bf p})=\left(\sqrt{m^2c^2+p^2},{\bf p}\right)$, $\theta(.)$ the Heaviside step function
and $\delta(.)$ the Dirac delta function.
It is trivial to verify that the latter integral transforms into the one-dimension integral
\begin{eqnarray}
I \propto \int d_{\kappa}x \, f(x) \ , \label{MIII10}
\end{eqnarray}
being $f(x)=4\pi x^2 F(x)$. Then the $\kappa$-integral is essentially the Lorentz invariant integral of special relativity.

\sect{The function $\exp_{\kappa}(x)$}

\subsection{Definition}

We recall that the ordinary exponential $f(x)=\exp(x)$ emerges as
solution both of the functional equation $f(x+y)=f(x)f(y)$ and of
the differential equation $(d/dx)f(x)=f(x)$. The question to
determine the solution of the  generalized equations
\begin{eqnarray}
f(x\stackrel{\kappa}{\oplus}y)=f(x)f(y) \ \ , \label{MIV1}
\end{eqnarray}
\begin{equation}
\frac{d\,f(x)}{d_{\kappa}x}=f(x) \ \ , \label{MIV2}
\end{equation}
reducing in the $\kappa \rightarrow 0$ limit to the ordinary
exponential, naturally arises. This solution is unique and
represents a one-parameter generalization of the ordinary
exponential.

{\it Solution of Eq. (\ref{MIV1}):} We write this equation
explicitly
\begin{eqnarray}
f\left(x\sqrt{1+\kappa^2y^2}+y\sqrt{1+\kappa^2x^2}\,\right)=f(x)f(y)
\ \ , \label{MIV3}
\end{eqnarray}
which after performing the change of variables $f(x)=\exp(g(\kappa
x))$, $z_1=\kappa x$, $z_2=\kappa y$ transforms as
\begin{eqnarray}
g\left(z_1\sqrt{1+z_2^2}+z_2\sqrt{1+z_1^2}\,\right)=g(z_1)+g(z_2) \
\ , \label{MIV4}
\end{eqnarray}
and admits as solution $g(x)=A\, {\rm arcsinh} x$. Then it results
that $f(x)=\exp(A\, {\rm arcsinh} \,\kappa x)$. The arbitrary
constant $A$ can be fixed through the condition $\lim_{\kappa
\rightarrow 0} f(x)=\exp(x)$, obtaining $A=1/\kappa$. Therefore
$f(x)$ assumes the form  $f(x)=\exp_{\kappa}\!\left(x \right)$ being
\begin{eqnarray}
\exp_{\kappa}(x)= \exp\left( \frac{1}{\kappa} \,{\rm arcsinh}\,
\kappa x \right) \ \ . \label{MIV5}
\end{eqnarray}

{\it Solution of Eq. (\ref{MIV2}):}

According to Eq. (\ref{MIV2}) the function $f(x)=\exp_{\kappa}(x)$
is defined as eigenfunction of $d/d_{\kappa}x$ i.e.
\begin{equation}
\frac{d\,\exp_{\kappa}(x)}{d_{\kappa}x}=\exp_{\kappa}(x) \ \ .
\label{MIV6} \nonumber
\end{equation}
After recalling that $d_{\kappa}x=d\{x\}$ with $\{x\}={\kappa}^{-1}
\,{\rm arcsinh}\, \kappa x $, Eq. (\ref{MIV6}) can be written in the
form
\begin{equation}
\frac{d\,\exp_{\kappa}(x)}{d\{x\}}=\exp_{\kappa}(x) \ \ .
\label{MIV7}
\end{equation}
The solution of the latter equation with the condition
$\exp_{\kappa}(0)=1$ follows immediately
\begin{eqnarray}
\exp_{\kappa}(x)=\exp \,(\{x\}) \ \ . \label{MIV8}
\end{eqnarray}
After taking into account that ${\rm arcsinh} \,x =\ln
(\sqrt{1+x^2}+x)$ we can write $\exp_{\kappa}(x)$ in the form
\begin{eqnarray}
\exp_{\kappa}(x)= \left(\sqrt{1+\kappa^{\,2} x^{\,2}}+\kappa
x\right)^{1/\kappa}\ \  , \label{MIV9}
\end{eqnarray}
which will used in the following.  We remark that $\exp_{\kappa}(x)$
given by Eq. (\ref{MIV9}), is solution both of the Eqs. (\ref{MIV1})
and (\ref{MIV2}) and therefore represents a generalization of the
ordinary exponential.

\subsection{Basic Properties}

From the definition (\ref{MIV9}) of $\exp_{\kappa}(x)$ follows that
\begin{eqnarray}
&&\exp_{\,0}(x)\equiv \lim_{\kappa \rightarrow
0}\exp_{\,\kappa}(x)=\exp (x) \ \ , \label{MIV10}
\\
&&\exp_{- \kappa}(x)=\exp_{\kappa}(x) \ \ .  \label{MIV11}
\end{eqnarray}
Like the ordinary exponential, $\exp_{\kappa}(x)$ has the
properties
\begin{eqnarray}
&&\exp_{\kappa}(x) \in C^{\infty}({\bf R}),
\label{MIV12}  \\
&&\frac{d}{d\,x}\, \exp_{\kappa}(x)>0,
\label{MIV13}  \\
&&\exp_{\kappa}(-\infty)=0^+,
\label{MIV14}  \\
&&\exp_{\kappa}(0)=1,
\label{MIV15}  \\
&&\exp_{\kappa}(+\infty)=+\infty,
\label{MIV16}  \\
&&\exp_{\kappa}(x)\exp_{\kappa}(-x)= 1 \ \ . \label{MIV17}
\end{eqnarray}

\begin{figure}[h]
\centerline{
\includegraphics[width=.8\columnwidth,angle=0]
{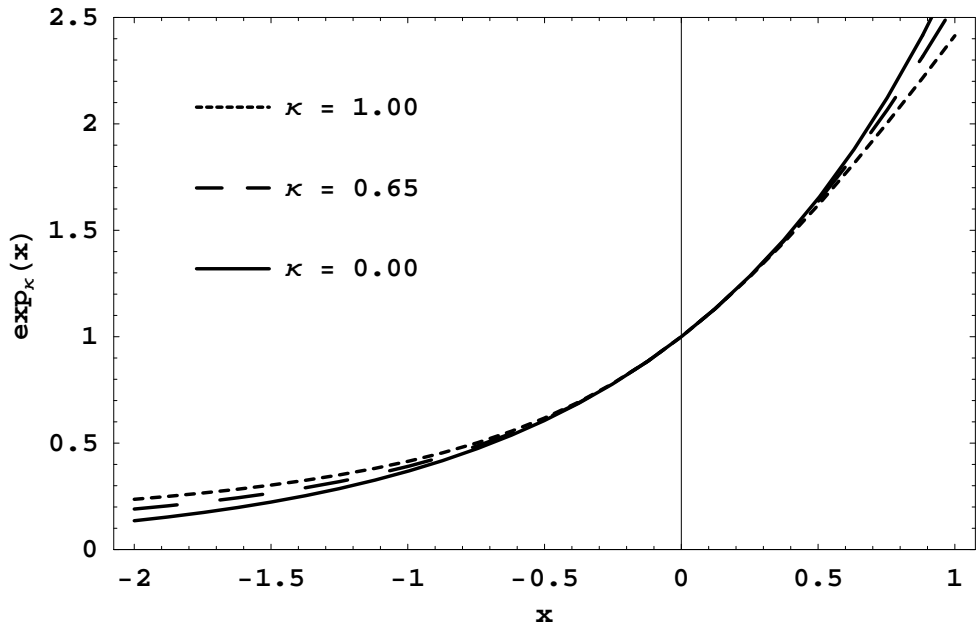}} \caption{Plot of  the function $\exp_{\kappa}(x)$ defined in Eq.(\ref{MIV9})
for three different values of the parameter of $\kappa$.
The continuous curve corresponding to $\kappa=0$ is the ordinary exponential function $\exp(x)$.}\label{FigA1}
\end{figure}

In Fig.\ref{FigA1} is plotted  the function $\exp_{\kappa}(x)$ defined in Eq.(\ref{MIV9})
for three different values of the parameter of $\kappa$. The continuous curve corresponding
to $\kappa=0$ is the ordinary exponential function $\exp(x)$.

The property (\ref{MIV17}) emerges as particular case of the more
general one
\begin{eqnarray}
\exp_{\kappa}(x)\exp_{\kappa}(y)=\exp_{\kappa}(x\stackrel{\kappa}{\oplus}y)\
\ . \label{MIV18}
\end{eqnarray}

Furthermore $\exp_{\kappa}(x)$ has the property
\begin{eqnarray}
&&\left (\exp_{\kappa}(x)\right )^{r} =\exp_{\kappa/r}(r x) \ \ ,
 \label{MIV19}
\end{eqnarray}
with $r\in {\bf R}$, which in the limit $\kappa \rightarrow 0$
reproduces one well known property of the ordinary exponential.

We remark the following convexity property
\begin{eqnarray}
\frac{d^2}{d\,x^2}\, \exp_{\kappa}(x)>0 \ \ ; \ \ x \in {\bf R} \ \
, \label{MIV20}
\end{eqnarray}
holding when $\kappa^2<1$.

Undoubtedly one of the more interesting properties of
$\exp_{\kappa}(x)$ is its power law asymptotic behavior
\begin{eqnarray}
\exp_{\kappa}(x) {\atop\stackrel{\textstyle\sim}{\scriptstyle
x\rightarrow \pm\infty}}\big|\,2\kappa x\big|^{\pm1/|\kappa|}
\ \ . \label{MIV21}
\end{eqnarray}

\subsection{Mellin transform}

Let us consider the incomplete Mellin transform of the $\exp_{\kappa}(-t)$
\begin{eqnarray}
{\cal M}_{\kappa}(r,x)&&=\int_{\, 0}^{x}\!t^{r-1} \,\exp_{\kappa}(-t)\,dt \ \ . \label{MIV22}
\end{eqnarray}
After performing the change of integration variable
$y=\left(\sqrt{1+\kappa^2t^2}-|\kappa| t \right) ^{2}$, and after
taking into account that
$t=\frac{1}{2|\kappa|}\left(y^{-1/2}-y^{1/2}\right)$ and
$\exp_{\kappa}(-t)=y^{1/2|\kappa|}$, the function ${\cal
M}_{\kappa}(r,x)$ can be written in the form
\begin{eqnarray}
{\cal M}_{\kappa}(r,x)&&= \frac{1}{2}\, |2\kappa|^{-r}\int_{X}^{1}\!
y^{\frac{1}{2|\kappa|}-\frac{r}{2}-1}
\left(1 -y\right)^{r-1}\,
\left(1+y \right)\, dy \ \ , \label{MIV23}
\end{eqnarray}
with
\begin{eqnarray}
X=\left(\sqrt{1+\kappa^2x^2}-|\kappa| x \right)
^{2} \ \ . \label{MIV24}
\end{eqnarray}

When $r$ is an integer greater that zero, ${\cal M}_{\kappa}(r,x)$
can be calculated analytically. For instance it results
\begin{eqnarray}
&&{\cal M}_{\kappa}(1,x)= \frac{1}{1-\kappa^2}-
\frac{\kappa^2x+\sqrt{1+\kappa^2x^2}}
{1-\kappa^2}\,\exp_{\kappa}(-x)\ \ , \label{MIV25} \\
&&{\cal M}_{\kappa}(2,x)= \frac{1}{1-4\kappa^2}-
\frac{1+2\kappa^2x^2+x\sqrt{1+\kappa^2x^2}}
{1-4\kappa^2}\,\exp_{\kappa}(-x) \ \ . \label{MIV26}
\end{eqnarray}

In general the function ${\cal M}_{\kappa}(r,x)$ can be written as
\begin{eqnarray}
{\cal M}_{\kappa}(r,x)=\frac{1}{2}\, |2\kappa|^{-r}\, \Big[I_1(r)-I_X(r)\Big] \ \ , \label{MIV27}
\end{eqnarray}
with
\begin{eqnarray}
I_X(r)=\int_{0}^{X}\!
y^{\frac{1}{2|\kappa|}-\frac{r}{2}-1}
\left(1 -y\right)^{r-1}\, dy +
\int_{0}^{X}\!
y^{\frac{1}{2|\kappa|}-\frac{r}{2}}
\left(1 -y\right)^{r-1}\, dy \ \ . \label{MIV28}
\end{eqnarray}
After recalling the definition of the Beta incomplete function $B_X
\left(s,r\right)=\int_{0}^{X}\! y^{s-1} \left(1 -y\right)^{r-1}\,
dy$, the integral $I_X(r)$ becomes
\begin{eqnarray}
I_X(r)=B_X \left(\frac{1}{2|\kappa|}-\frac{r}{2}\,,\,r\right)
+B_X \left(\frac{1}{2|\kappa|}-\frac{r}{2}+1\,,\,r\right) \ \ . \label{MIV29}
\end{eqnarray}
The function $I_1(r)$ can be expressed in term of the Beta functions
$B\left(s,r\right)=\int_{0}^{1}\! y^{s-1} \left(1 -y\right)^{r-1}\,
dy$, and then in terms of Gamma functions, obtaining
\begin{eqnarray}
I_1(r)=\frac{2\, \Gamma\left(r\right)}{1+|\kappa| r}\,
\frac{\Gamma\left(\frac{1}{2|\kappa|}-\frac{r}{2}\right)
}{\Gamma\left(\frac{1}{2|\kappa|}+\frac{r}{2}\right)} \ \ . \label{MIV30}
\end{eqnarray}
Finally the incomplete Mellin transform ${\cal M}_{\kappa}(r,x)$ of $\exp_{\kappa}(-t)$ assumes the form
\begin{eqnarray}
\nonumber {\cal M}_{\kappa}(r,x)&&=\frac{|2\kappa|^{-r}}{1+|\kappa| r}\,
\frac{\Gamma\left(\frac{1}{2|\kappa|}-\frac{r}{2}\right)
}{\Gamma\left(\frac{1}{2|\kappa|}+\frac{r}{2}\right)}\, \Gamma\left(r\right)
\\ \nonumber &&-\frac{1}{2}\,|2\kappa|^{-r}\,B_X \left(\frac{1}{2|\kappa|}-\frac{r}{2}\,,\,r\right)
\\&&-\frac{1}{2}\,|2\kappa|^{-r}\,B_X \left(\frac{1}{2|\kappa|}-\frac{r}{2}+1\,,\,r\right)
\ \ . \label{MIV31}
\end{eqnarray}

The Mellin transform of $\exp_{\kappa}(-t)$, namely
\begin{eqnarray}
{\cal M}_{\kappa}(r)=\int_{\, 0}^{\infty}\!t^{r-1} \,\exp_{\kappa}(-t)\,dt
\ \ , \label{MIV32}
\end{eqnarray}
can be calculated from Eq. (\ref{MIV31}) by posing  $x=\infty$. The
explicit expression of ${\cal M}_{\kappa}(r)$ holding for $0<
r<1/|\kappa|$ is given by

\begin{eqnarray}
{\cal M}_{\kappa}(r)
=\frac{|2\kappa|^{-r}}{1+|\kappa| r}\, \frac{\Gamma\left(\frac{1}{2|\kappa|}-\frac{r}{2}\right)
}{\Gamma\left(\frac{1}{2|\kappa|}+\frac{r}{2}\right)}\, \Gamma \left(r\right)
\ \ . \label{MIV33}
\end{eqnarray}

From the latter relationship one can verify easily the property
\begin{eqnarray}
{\cal M}_{\kappa}(r+2)=\frac{r(r+1)}{1-\kappa^2\,(r+2)^2}\,{\cal M}_{\kappa}(r)
\ \ . \label{MIV34}
\end{eqnarray}

\subsection{Taylor expansion}

The Taylor expansion of $\exp_{\kappa}(x)$ given in \cite{PRE2002} can be written also in the following form
\begin{equation}
\exp_{\kappa}(x) = \sum_{n=0}^{\infty}\frac{x^n}{n!_{\kappa}} \ \ \
; \ \ \ \kappa^2 x^2 < 1 \ \ ,
 \label{MIV35}
\end{equation}
where the symbol $n!_{\kappa}$, representing the $\kappa$-generalization of the
ordinary factorial $n!$, recovered for $\kappa=0$, is given by
\begin{equation}
n!_{\kappa}=\frac{n!}{\xi_n(\kappa)} \ \ , \label{MIV36}
\end{equation}
and the polynomials $\xi_n(\kappa)$ are defined as
\begin{eqnarray}
&&\xi_0(\kappa)=\xi_1(\kappa)=1 \ \ ,  \label{MIV37}  \\
&&\xi_n(\kappa)=\prod_{j=1}^{n-1}\Big[1-(2j-n)\,\kappa\Big] \ \ ;
\ \ n>1 \ \ . \label{MIV38}
\end{eqnarray}

The polynomials $\xi_n(\kappa)$, for $n>1$, when $n$ is odd, are of degree $n-1$,
with respect the variable $\kappa$, while when $n$ is even the
degree of $\xi_n(\kappa)$ is $n-2$. The degree of $\xi_n(\kappa)$ is
always an even number and it results
\begin{eqnarray}
&&\xi_{2m}(\kappa)=\prod_{j=0}^{m-1}\Big[1-(2j)^2\kappa^2\Big] \ \
; \ \ m>0 \ \ , \label{MIV39} \\
&&\xi_{2m+1}(\kappa)= \prod_{j=0}^{m-1}\Big[1-(2j+1)^2\kappa^2\Big]
\ \  ; \ \ m>0 \ \ .  \label{MIV40}
\end{eqnarray}

The polynomials $\xi_n(\kappa)$ can be generated by the following
simple recursive formula
\begin{eqnarray}
&&\xi_0(\kappa)=\xi_1(\kappa)=1 \ \ ,  \label{MIV41}  \\
&&\xi_{n+2}(\kappa)=(1-n^2\kappa^2)\,\xi_n(\kappa) \ \ ; \ \ n\geq 0
\ \ . \label{MIV42}
\end{eqnarray}

The first nine polynomials reads as
\begin{eqnarray}
&&\xi_0(\kappa)=\xi_1(\kappa)=\xi_2(\kappa)=1 \ \ ,  \label{MIV43}  \\
&&\xi_{3}(\kappa)=1-\kappa^2 \ \ , \label{MIV44} \\
&&\xi_{4}(\kappa)=1-4\kappa^2 \ \ , \label{MIV45} \\
&&\xi_{5}(\kappa)=(1-\kappa^2)(1-9\kappa^2) \ \ , \label{MIV46} \\
&&\xi_{6}(\kappa)=(1-4\kappa^2)(1-16\kappa^2) \ \ , \label{MIV47} \\
&&\xi_{7}(\kappa)=(1-\kappa^2)(1-9\kappa^2)(1-25\kappa^2) \ \ , \label{MIV48} \\
&&\xi_{8}(\kappa)=(1-4\kappa^2)(1-16\kappa^2)(1-36\kappa^2) \ \ .
\label{MIV49}
\end{eqnarray}

After noting that for a given value of $\kappa$ the maximum natural
number $N$ satisfying the condition $N<2+1/|\kappa|$ is defined
univocally, we can verify easily that for $n=0,1,2,...,N$ it results
$\xi_n(\kappa)>0$ and then $n!_{\kappa}>0$. For $n> N$ the sign of
$\xi_n(\kappa)$ and then of $n!_{\kappa}$ alternates with
periodicity $--++--++...$

From Eqs. (\ref{MIV36}) and (\ref{MIV42}) follows the recursive
formula
\begin{equation}
(n+2)!_{\,\kappa}=\frac{(n+1)(n+2)}{1-n^2\kappa^2}\,n!_{\kappa}\ \ .
\label{MIV50}
\end{equation}
By direct comparison of Eq. (\ref{MIV34}) and (\ref{MIV50}) we
obtain the relationship
\begin{eqnarray}
n!_{\kappa}=\left(1-\kappa^2n^2\,\right)n\int_{\, 0}^{\infty}\!t^{n-1} \,\exp_{\kappa}(-t)\,dt
 \ \ . \label{MIV51}
\end{eqnarray}

It is remarkable that the first three terms in the Taylor expansion
of $\exp_{\kappa}(x)$ are the same of the ordinary exponential
\begin{equation}
\exp_{\kappa}(x) = 1+ x + \frac{x^2}{2} +
(1-\kappa^2)\,\frac{x^3}{3!}+...  \ \ . \label{MIV52}
\end{equation}

\subsection{The function $\Gamma_{\kappa}(x)$}

The $\Gamma_{\kappa}(n)$ with $n$ integer is defined through
\begin{eqnarray}
\Gamma_{\kappa}(n)= (n-1)!_{\kappa} \ \ ,  \ \ \ \ \label{MIV53}
\end{eqnarray}
and represents a generalization of the Euler $\Gamma(n)$ function.
In particular we have $\Gamma_{\kappa}(1)=\Gamma_{\kappa}(2)=1$ and
$\Gamma_{\kappa}(3)=2$. This definition and the relationship
(\ref{MIV51}) suggests the following one parameter generalization of
the Euler $\Gamma(x)$ function i.e. $\Gamma_{\kappa}(x)$, given by
\begin{eqnarray}
\Gamma_{\kappa}(x)=\left[1-\kappa^2(x-1)^2\right](x-1)\int_{\, 0}^{\infty}\!t^{x-2} \,
\exp_{\kappa}(-t)\,dt
 \ \ . \label{MIV54}
\end{eqnarray}
The explicit expression of $\Gamma_{\kappa}(x)$ in terms of the ordinary $\Gamma(x)$ is given by
\begin{eqnarray}
\Gamma_{\kappa}(x)=\frac{1-|\kappa|(x-1)}{|2\kappa|^{x-1}}\,\,\,
\frac{\Gamma\left(\frac{1}{|2\kappa|}-\frac{x-1}{2}
\right)}{\Gamma\left(\frac{1}{|2\kappa|}+\frac{x-1}{2}
\right)}\,\,\Gamma\left(x\right)\ , \ \ \ \ \label{MIV55}
\end{eqnarray}
and can be used as definition of  $\Gamma_{\kappa}(x)$ when $x$ is a complex variable.
Clearly in the $\kappa \rightarrow 0$ limit it results $\Gamma_{0}(x)=\Gamma(x)$.
An expression of $\Gamma_{\kappa}(x)$ in terms of the Beta function is the following
\begin{eqnarray}
\Gamma_{\kappa}(x)=\frac{1-|\kappa|(x-1)}{|2\kappa|^{x-1}}\,\,(x-1)\,
B\left(\frac{1}{|2\kappa|}-\frac{x-1}{2} \, , \,  x-1 \right)\  . \
\ \ \ \label{MIV56}
\end{eqnarray}

From Eqs. (\ref{MIV34}) and (\ref{MIV54}) follows the property
\begin{eqnarray}
\Gamma_{\kappa}(x+2)=\frac{x\,(x+1)}{1-\kappa^2(x-1)^2}\,
\Gamma_{\kappa}(x) \ \ . \label{MIV57}
\end{eqnarray}

The Taylor expansion of $\exp_{\kappa}(x)$ can be written also in the form
\begin{equation}
\exp_{\kappa}(x) = \sum_{n=0}^{\infty}\,\frac{x^n}{\Gamma_{\kappa}(n+1)} \ \ \
; \ \ \ \kappa^2 x^2 < 1 \ \ .
 \label{MIV58}
\end{equation}

\begin{figure}[h]
\centerline{
\includegraphics[width=.8\columnwidth,angle=0]
{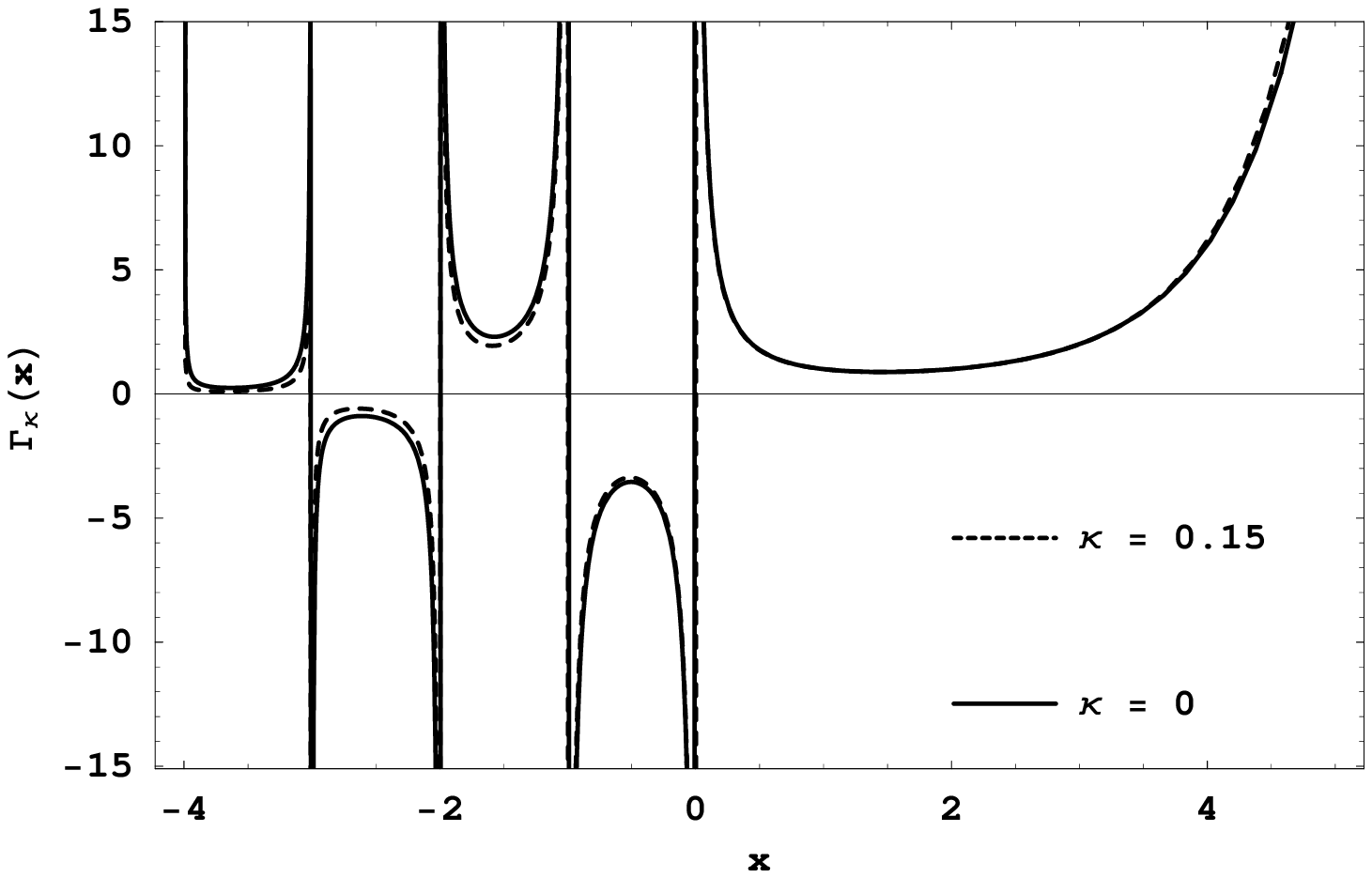}} \caption{Plot of  the function $\Gamma_{\kappa}(x)$ defined in Eq.(\ref{MIV55})
in the range  $-4<x<4$  for $\kappa=0$ and $\kappa=0.15$. The continuous curve corresponding to $\kappa=0$
is the ordinary Gamma function $\Gamma(x)$.}\label{FigA2}
\end{figure}

\begin{figure}[h]
\centerline{
\includegraphics[width=.8\columnwidth,angle=0]
{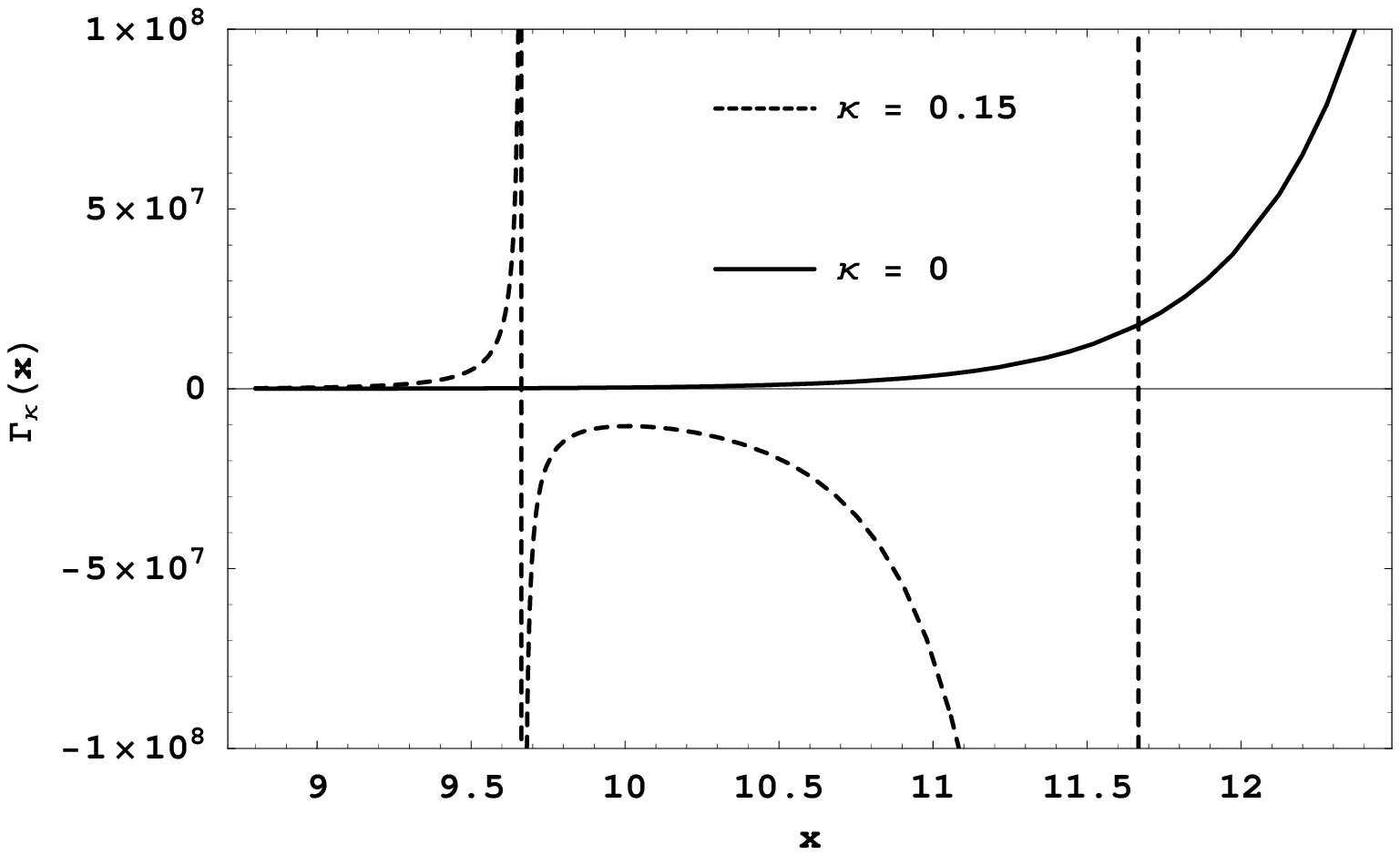}} \caption{Plot of  the function $\Gamma_{\kappa}(x)$ defined in Eq.(\ref{MIV55})
in the range  $9<x<12$  for $\kappa=0$ and $\kappa=0.15$. The continuous curve corresponding to $\kappa=0$
is the ordinary Gamma function $\Gamma(x)$.}\label{FigA3}
\end{figure}

In Fig.\ref{FigA2} and Fig.\ref{FigA3} is plotted the function
$\Gamma_{\kappa}(x)$ defined in Eq.(\ref{MIV55}) in the ranges
$-4<x<4$ and $9<x<12$ respectively, for $\kappa=0$ and
$\kappa=0.15$. The continuous curve corresponding to $\kappa=0$ is
the ordinary Gamma function $\Gamma(x)$.

The incomplete $\gamma_{\kappa}(r,x)$ and $\Gamma_{\kappa}(r,x)$ are defined as
\begin{eqnarray}
\gamma_{\kappa}(r,x)=\left[1-\kappa^2(r-1)^2\right](r-1)\int_{\, 0}^{x}\!t^{r-2} \,
\exp_{\kappa}(-t)\,dt
 \ \ , \\  \label{MIV59}
\Gamma_{\kappa}(r,x)=\left[1-\kappa^2(r-1)^2\right](r-1)\int_{\, x}^{\infty}\!t^{r-2} \,
\exp_{\kappa}(-t)\,dt
 \ \ ,  \label{MIV60}
\end{eqnarray}
and hold the following relationships
\begin{eqnarray}
&&\gamma_{\kappa}(r,x)+\Gamma_{\kappa}(r,x)=\Gamma_{\kappa}(r)
 \ \   \label{MIV61} \ \ , \\
&&\gamma_{\kappa}(r,\infty)=\Gamma_{\kappa}(r)
 \ \   \label{MIV62} \ \ , \\
&&\gamma_{\kappa}(r,x)=(r-1)\!\left[1-\kappa^2(r-1)^2\,\right]{\cal M}_{\kappa}(r-1,x)
 \ \   \label{MIV63} \ \ .
\end{eqnarray}

\subsection{Expansion in ordinary exponentials}

Starting from the expression (\ref{MIV5}) $\exp_{\kappa}(x)$ and the Taylor expansion of the function
${\rm arcsinh}(x)$  we obtain
\begin{eqnarray}
\exp_{\kappa}(x)= \exp \Big(\sum_{n=0}^{\infty}
c_n\,\kappa^{2n}x^{2n+1}\Big) \ \ \ \ , \ \ \  \kappa^2x^2\leq 1 \ \
\ ,  \label{MIV64}
\end{eqnarray}
with
\begin{eqnarray}
c_n=\frac{(\!-1)^n\,(2n)!}{(2n\!+\!1)\,2^{2n}\,(n!)^2} \ \ .
\label{MIV65}
\end{eqnarray}
Exploiting this relationship, we can write $\exp_{\kappa}(x)$ as
an infinite product of ordinary exponentials
\begin{eqnarray}
\exp_{\kappa}(x)=\prod_{n=0}^{\infty} \exp
\Big(c_n\,\kappa^{2n}x^{2n+1}\Big) \ \ . \label{MIV66}
\end{eqnarray}

On the other hand $\exp_{\kappa}(x)$ can be viewed as a continuous
linear combination of an infinity of standard exponentials. Namely
for ${\rm Re}\,s \geq 0$, it holds the following Laplace transform

\begin{equation}
\exp_{\kappa}(-s) = \int_0^{\infty}\frac{1}{\kappa
x}\,\,J_{_{\scriptstyle \!1/\kappa}}\Big(\frac{x}{\kappa}\Big) \exp
\,(-sx)\, dx \ \ ,  \label{MIV67}
\end{equation}
being $J_{\nu}(x)$ the Bessel function.

\subsection{$\kappa$-Laplace transform}

The following $\kappa$-Laplace transform emerges naturally
\begin{eqnarray}
F_{\kappa}(s)={\cal L}_{\kappa}\{f(t)\}(s)=\int_{\, 0}^{\infty}\!f(t) \,[\exp_{\kappa}(-t)]^s\,dt  \ \ , \label{MIV68}
\end{eqnarray}
as a generalization of the ordinary Laplace transform. The inverse $\kappa$-Laplace transform is given by
\begin{eqnarray}
f(t)={\cal L}^{-1}_{\kappa}\{F_{\kappa}(s)\}(t)={\frac{1}{2\pi i}\int_{c-i \infty}^{c+i \infty}\!F_{\kappa}(s) \,\frac{[\exp_{\kappa}(t)]^s}{\sqrt{1+\kappa^2t^2}}\,ds}  \ \ . \label{MIV69}
\end{eqnarray}

In ref. \cite{DeossaCasas} the mathematical properties of the $\kappa$-Laplace transform have been investigated systematically. In table I are reported the main properties of the $\kappa$-Laplace transform which in the $\kappa \rightarrow 0$ limit reduce to the corresponding ordinary Laplace transform properties.

\begin{table*}
\caption{\label{tab:table1} Properties of the $\kappa$-Laplace transform.}

\begin{ruledtabular}
\begin{tabular}{cc}
\\
$f(t)$
&
$F_{\kappa}(s)$
\\ \\

\hline \hline
\\
$a\, f (t)+ b\, g (t)$
&
$a\, F_{\kappa} (s)+ b\, G_{\kappa} (s) $
\\ \\
\hline
\\
$ f(at)$
&
$\frac{1}{a}\, F_{\kappa / a} (\frac{s}{a}) $
\\ \\
\hline
\\
$f(t)\, [\exp_{\kappa}(-t)]^{a} $
&
$F_{\kappa}(s-a) $
\\ \\
\hline
\\
$\frac{d\, f(t)}{dt} $
&
$s\, {\cal L}_{\kappa}\left \{\frac{f(t)}{\sqrt{1+\kappa^2 t^2}}\right \}(s)-f(0)  $
\\ \\
\hline
\\
$\frac{d}{dt} \,  \sqrt{1+\kappa^2 t^2} \, f(t) $
&
$ s \, F_{\kappa} (s) -f(0)$
\\ \\
\hline
\\
$\frac{1}{\sqrt{1+\kappa^2 t^2}}\, \int_0^t f(w)dw$
&
$\frac{1}{s} \, F_{\kappa} (s) $
\\ \\
\hline\\
$f(t)\, [\ln (\exp_{\kappa}(t))]^n $
&
$(-1)^n \frac{d^{n} F_{\kappa}(s)}{ds^n}$
\\ \\
\hline
\\
$ f(t) \,[\ln (\exp_{\kappa}(t))]^{-n}$
&
$ \int_s^{+\infty}dw_{n} \int_{w_n}^{+\infty}dw_{n-1}...\int_{w_3}^{+\infty}dw_{2}\int_{w_2}^{+\infty}dw_{1} \,F_{\kappa}(w_1)$
\\
\\
\end{tabular}
\end{ruledtabular}
\end{table*}

Furthermore it holds the initial value theorem
\begin{eqnarray}
\lim_{t\rightarrow 0}f(t)=\lim_{s\rightarrow \infty}s F_{\kappa}(s) \ \ , \label{MIV70}
\end{eqnarray}
and the final value theorem
\begin{eqnarray}
\lim_{t\rightarrow \infty}|\kappa|t f(t)=\lim_{s\rightarrow 0}s F_{\kappa}(s) \ \ . \label{MIV71}
\end{eqnarray}

The $\kappa$-convolution of two functions $f \stackrel{\kappa}{*} g =(f \stackrel{\kappa}{*} g )(t)$, is defined as
\begin{eqnarray}
f \stackrel{\kappa}{*} g =\int_0^t f(t\stackrel{\kappa}{\ominus}\tau) \, g(\tau) \, \frac{1-\kappa^2\tau (t-\tau)}{\sqrt{1+\kappa^2 \tau^2}} \, d \tau \ \ . \label{MIV72}
\end{eqnarray}
and has the following properties
\begin{eqnarray}
&&f \stackrel{\kappa}{*} (ag+bh)=a (f \stackrel{\kappa}{*} g)+b (f \stackrel{\kappa}{*} h) \ , \label{MIV73} \\
&&f \stackrel{\kappa}{*} g=g \stackrel{\kappa}{*} f  \ ,  \label{MIV74}\\
&&f \stackrel{\kappa}{*} (g\stackrel{\kappa}{*}h)=(f \stackrel{\kappa}{*} g) \stackrel{\kappa}{*} h \label{MIV75} \ .
\end{eqnarray}
It holds the following $\kappa$-convolution theorem
\begin{eqnarray}
{\cal L}_{\kappa}\{f \stackrel{\kappa}{*} g\}= {\cal L}_{\kappa}\{f\} \, {\cal L}_{\kappa}\{g\} \ \ . \label{MIV68}
\end{eqnarray}

In table II are reported the $\kappa$-Laplace transforms for the delta function, for the unit function and for the power function. We note that ${\kappa}$-Laplace transform of the power function $f(t)=t^{\nu -1}$ involves the $\kappa$-generalized Gamma function. All the $\kappa$-Laplace transforms in the $\kappa \rightarrow 0$ limit reduce to the corresponding ordinary Laplace transforms.

\begin{table*}
\caption{\label{tab:table1} The $\kappa$-Laplace transform of the Dirac delta-function, of the Heaviside unit function and of the power function.}

\begin{ruledtabular}
\begin{tabular}{cc}
\\
$f(t)$
&
$F_{\kappa}(s)$
\\ \\

\hline \hline
\\
$\delta (t-\tau)$
&
$[\exp_{\kappa}(-\tau)]^s$
\\ \\
\hline
\\
$u(t-\tau)$
&
$\frac{s\sqrt{1+\kappa^2 \tau^2}+\kappa^2 \tau}{s^2-\kappa^2}\, [\exp_{\kappa}(-\tau)]^{s} $
\\ \\
\hline
\\
$t^{\nu-1}$
&
$\frac{s^2}{s^2-\kappa^2\nu^2}\,\frac{\Gamma_{\frac{\kappa}{s}}(\nu+1)}{\nu\, s^{\nu}}=\frac{s}{s+|\kappa|\nu}\, \frac{\Gamma (\nu)}{|2\kappa|^{\nu}}\, \frac{\Gamma\left( \frac{s}{|2\kappa|} - \frac{\nu}{2} \right )} {\Gamma\left( \frac{s}{|2\kappa|} + \frac{\nu}{2} \right )}$
\\ \\
\hline
\\
$t^{2m-1}, \ \ m \in Z^+ $
&
$\frac{(2m-1)!}{\prod_{j=1}^{m}\left[s^2-(2j)^2\kappa^2\right] }$
\\ \\
\hline
\\
$t^{2m}, \ \ m \in Z^+ $
&
$\frac{(2m)!\, s}{\prod_{j=1}^{m+1}\left[s^2-(2j-1)^2\kappa^2\right] }$
\\ \\
\end{tabular}
\end{ruledtabular}
\end{table*}

\sect{The function $\ln_{\kappa}(x)$}

\subsection{Definition and basic properties}

The function $\ln_{\kappa}(x)$ is defined as the inverse function
of $\exp_{\kappa}(x)$, namely
\begin{eqnarray}
\ln_{\kappa}(\exp_{\kappa}x)=\exp_{\kappa}(\ln_{\kappa}x)=x \ ,
\label{MV1}
\end{eqnarray}
and is given by
\begin{eqnarray}
\ln_{\kappa}(x)= [\ln x] \ \ , \label{MV2}
\end{eqnarray}
and then
\begin{eqnarray}
\ln_{\kappa}(x) = \frac{1}{\kappa}\,\sinh \, (\kappa \ln x)\ \ ,
\label{MV3}
\end{eqnarray}
or more properly
\begin{eqnarray}
\ln_{\kappa}(x)= \frac{x^{\kappa}-x^{-\kappa}}{2\kappa} \ \ .
\label{MV4}
\end{eqnarray}

It results that
\begin{eqnarray}
&&\ln_{0}(x)\equiv \lim_{\kappa \rightarrow 0}\ln_{\kappa}(x)=\ln
(x) \ \ , \label{MV5}
\\
&&\ln_{- \kappa}(x)=\ln_{\kappa}(x) \ \ . \label{MV6}
\end{eqnarray}

The function $\ln_{\kappa}(x)$, just as the ordinary logarithm, has
the properties
\begin{eqnarray}
&&\ln_{\kappa}(x) \in C^{\infty}({\bf R}^+),
\label{MV7} \\
&&\frac{d}{d\,x}\, \ln_{\kappa}(x)>0,
\label{MV8}  \\
&&\ln_{\kappa}(0^+)=-\infty,
\label{MV9}  \\
&&\ln_{\kappa}(1)=0,
\label{MV10}  \\
&&\ln_{\kappa}(+\infty)=+\infty,
\label{MV11}  \\
&&\ln_{\kappa}(1/x)=-\ln_{\kappa}(x) \ \ . \label{MV12}
\end{eqnarray}

\begin{figure}[h]
\centerline{
\includegraphics[width=.8\columnwidth,angle=0]
{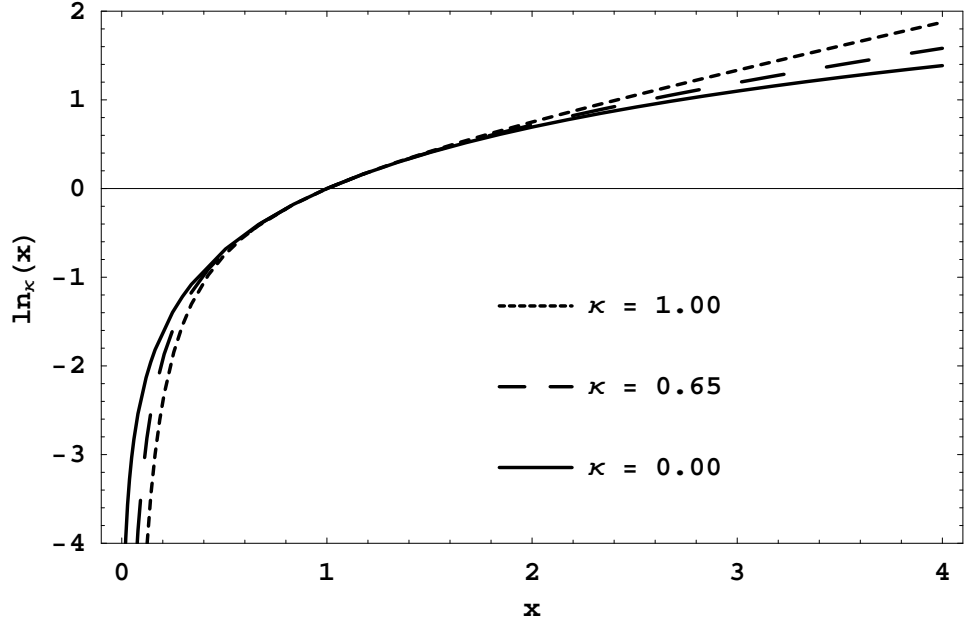}} \caption{Plot of  the function $\ln_{\kappa}(x)$ defined by Eq.(\ref{MV4})
for three different values of the parameter of $\kappa$.
The continuous curve corresponding to $\kappa=0$ is the ordinary logarithm function $\ln(x)$.}\label{FigA4}
\end{figure}

In Fig.\ref{FigA4} is plotted  the function $\ln_{\kappa}(x)$ defined in Eq.(\ref{MV4})
for three different values of the parameter of $\kappa$. The continuous curve corresponding
to $\kappa=0$ is the ordinary logarithm function $\ln(x)$.

Furthermore $\ln_{\kappa}(x)$ has the two properties
\begin{eqnarray}
&&\ln_{\kappa}(x^{r}) =r \ln_{r \kappa}(x) \ \ ,  \label{MV13}
\\
&&\ln_{\kappa}(x \,y) =\ln_{\kappa}(x)
\oplus\!\!\!\!\!^{^{\scriptstyle \kappa}}\,\,\ln_{\kappa}(y)\ \ ,
 \label{MV14}
\end{eqnarray}
with $r\in {\bf R}$. Note that the property (\ref{MV12}) follows as
particular case of the property (\ref{MV14}).

We remark the following concavity properties
\begin{eqnarray}
&&\frac{d^2}{d\,x^2}\, \ln_{\kappa}(x)<0
 \ \ , \label{MV15}   \\
&&\frac{d^2}{d\,x^2}\,x\, \ln_{\kappa}(x)<0  \ \ . \label{MV16}
\end{eqnarray}

A very interesting property of this function is its power law
asymptotic behavior
\begin{eqnarray}
&&\ln_{\kappa}(x) {\atop\stackrel{\textstyle\sim}{\scriptstyle
x\rightarrow
0^+}}-\frac{1}{2\,|\kappa|}\,x^{-|\kappa|} \ \ , \label{MV17}  \\
&&\ln_{\kappa}(x) {\atop\stackrel{\textstyle\sim}{\scriptstyle
x\rightarrow +\infty}}\,\,\frac{1}{2\,|\kappa|}\,x^{|\kappa|} \ \
.   \label{MV18}
\end{eqnarray}

After recalling the integral representation of the ordinary logarithm
\begin{equation}
\ln(x) = \frac{1}{2}\,\int_{1/x}^{x}\frac{1}{t} \,dt \ \ ,
\nonumber  \label{MV19}
\end{equation}
one can verify that the latter relationship can be generalized
easily in order to obtain $\ln_{\kappa}(x)$,  by replacing the
integrand function $y_0(t)=t^{-1}$ by the new function
$y_{\kappa}(t)=t^{-1-\kappa}$, namely
\begin{equation}
\ln_{\kappa}(x) = \frac{1}{2}\,\int_{1/x}^{x}\,\frac{1}{t^{1+\kappa}} \,dt \ \ .
\nonumber  \label{MV20}
\end{equation}

\subsection{Taylor expansion}

The Taylor expansion of $\ln_{\kappa}(1+x)$ converges if $-1<x\leq
1$, and assumes the form
\begin{equation}
\ln_{\kappa}(1+x) = \sum_{n=1}^{\infty} b_n(\kappa)\,(-1)^{n-1}
\,\frac{x^n}{n} \ \ , \label{MV21} \nonumber
\end{equation}
with $b_1(\kappa)= 1$, while for $n>1$, $b_n(\kappa)$ is given by
\begin{eqnarray}
b_{n}(\kappa)&&\,\,=\frac{1}{2}\Big(1-\kappa\Big)\Big(1-\frac{\kappa}{2}\Big)...
\Big(1-\frac{\kappa}{n-1}\Big) \nonumber
\\ &&\,\,+\,\frac{1}{2}\Big(1+\kappa\Big)\Big(1+\frac{\kappa}{2}\Big)...
\Big(1+\frac{\kappa}{n-1}\Big) \label{MV22}  \ \ .
\end{eqnarray}
It results $b_n(0)=1$ and $b_n(-\kappa)=b_n(\kappa)$. The first
terms of the expansion are
\begin{equation}
\ln_{\kappa}(1+x) = x - \frac{x^2}{2} +
\left(1+\frac{\kappa^2}{2}\right)\frac{x^3}{3} - ... \ \ .
\nonumber  \label{MV23}
\end{equation}

\subsection{The function $\Gamma_{\kappa}(x)$}

The following integral is useful

\begin{eqnarray}
&&\int_{\, 0}^{1}\left(\ln_{\kappa}\frac{1}{x}\right)^{r-1}\,dx \nonumber \\
&&=\frac{|2\kappa|^{1-r}}{1+(r-1)|\kappa|}\,\,\,\frac{\Gamma\left(\frac{1}{|2\kappa|}-\frac{r-1}{2}
\right)}{\Gamma\left(\frac{1}{|2\kappa|}+\frac{r-1}{2} \right)}
\,\,\Gamma\left(r\right) \ . \ \ \ \ \ \ \ \label{MV24}
\end{eqnarray}

Starting from the definition of the generalized Euler gamma
function, i.e. $\Gamma_{\kappa}(x)$ given in the previous section,
we can write it also in the following alternative but equivalent
form
\begin{eqnarray}
\Gamma_{\kappa}(x)=\left[1-\kappa^2(x-1)^2\right] \int_{\, 0}^{1}
\left(\ln_{\kappa}\frac{1}{t}\right)^{x-1} dt
 \ \ . \label{MV25}
\end{eqnarray}
An expression of $\Gamma_{\kappa}(x)$, where the parameter $\kappa$
enters exclusively through the function $\ln_{\kappa}(.)$, follows
easily
\begin{eqnarray}
\Gamma_{\kappa}(x)=(x-1)\,\,\frac{\int_{\, 0}^{1}
\left(\ln_{\kappa}\frac{1}{t}\right)^{x-1} dt}{\int_{\, 0}^{1}
\ln_{\kappa}\left(\frac{1}{t}\right)^{x-1} dt}
 \ \ . \label{MV26}
\end{eqnarray}

From the latter relationships follows that $n!_{\kappa}$ is given by

\begin{eqnarray}
n!_{\kappa}=\left(1-\kappa^2n^2\,\right) \int_{\, 0}^{1}
\left(\ln_{\kappa}\frac{1}{t}\right)^{n} dt=n\,\,\frac{\int_{\, 0}^{1}
\,\left(\ln_{\kappa}\frac{1}{t}\right)^n dt}{\int_{\, 0}^{1}\,
\ln_{\kappa}\left(\frac{1}{t}\right)^n dt}
 \ \ . \label{MV27}
\end{eqnarray}

\subsection{$\ln_{\kappa}(x)$ as solution of a functional equation}

The logarithm $y(x)=\ln(x)$ is the only existing function, unless a
multiplicative constant, which results to be solution of the
function equation $y(x_1x_2)=y(x_1)+y(x_2)$. Let us consider now the
generalization of this equation,  obtained by substituting the
ordinary sum by the generalized sum
\begin{eqnarray}
y(x_1x_2)= y(x_1)\stackrel{\kappa}{\oplus} y(x_2) \ \ . \label{MV28}
\end{eqnarray}
We proceed by solving this equation, which assumes the explicit
form
\begin{eqnarray}
y(x_1x_2)&=& y(x_1)\,\sqrt{1+\kappa^2\,y(x_2)\,^2} \nonumber \\
&+& y(x_2)\,\sqrt{1+\kappa^2\,y(x_1)\,^2} \ \ \ . \label{MV29}
\end{eqnarray}
After performing the substitution $y(x)=\kappa^{-1} \sinh \kappa
g(x)$ we obtain that the auxiliary function $g(x)$ obeys the
equation $g(x_1x_2)=g(x_1)+g(x_2)$, and then is given by $g(x)=A\ln
x$. The unknown function $y(x)$ becomes $y(x)=\kappa^{-1}\sinh (
\kappa \ln x)$ where we have set $A=1$ in order to recover, in the
limit $\kappa\rightarrow 0$, the classical solution $y(x)=\ln(x)$.
Then we can conclude that the solution of Eq. ($\ref{MV28}$) is
given by
\begin{eqnarray}
y(x)=\ln_{\kappa}(x) \ \ . \label{MV30}
\end{eqnarray}

\subsection{$\ln_{\kappa}(x)$ as solution of a differential-functional equation}

The following first order differential-functional equation emerges
in statistical mechanics within the context of the maximum entropy
principle
\begin{eqnarray}
&&\frac{d}{dx}\,[\,x\,\,f\,(x)\,]=\frac{1}{\gamma} \,f\left(\epsilon x\right) \
\ , \label{MV31}
\\&&f(1)=0 \ \ , \label{MV32}
\\ &&f'(1)=1 \ \ , \label{MV33}
\\ &&f\left(1/x\right)=-f\left(x\right). \label{MV34}
\end{eqnarray}
The latter problem admits two solutions \cite{PRE2002,PLA2011}. The first is given by
$f(x)=\ln (x)$ and $\gamma=1$, $\epsilon=e$. The second solution is
given by
\begin{eqnarray}
f(x)= \ln_{\kappa}(x) \ , \ \ \ \label{MV35}
\end{eqnarray}
and
\begin{eqnarray}
&&\gamma=\frac{1}{\sqrt{1-\kappa^2}} \ , \label{MV36} \\
&&\epsilon=\left(\frac{1+\kappa}{1-\kappa}\right)^{1/2\kappa} \ .
\label{MV37}
\end{eqnarray}

The constant $\gamma$ is the Lorentz factor corresponding to the reference velocity $v_*$ while the constant $\epsilon=\exp_{\kappa}(\gamma)$, represents the $\kappa$-generalization of the Napier number $e$.

\subsection{The Entropy}

A physically meaningful link between the functions $\ln_{\kappa}(x)$ and $\exp_{\kappa}(x)$ is given by a variational principle.
It holds the following theorem:

\noindent {\bf Theorem:} Let be $g(x)$ an arbitrary real function and $y(x)$
a probability distribution function of the variable $x\in A$.
The solution of the variational equation
\begin{eqnarray}
\frac{\delta}{\delta y(x)}\left[-\int_Adx \,\,y(x)\ln_{\kappa}y(x)
+\int_Adx \,\,y(x)\,g(x) \right]=0 \ ,
\label{MV38}
\end{eqnarray}
is unique and is given by
\begin{eqnarray}
y(x)=\frac{1}{\epsilon} \, \exp_{\kappa}\!\big(\gamma\, g(x)\big) \ ,
\label{MV39}
\end{eqnarray}
being the constants $\gamma$ and $\epsilon$  defined by Eqs. (\ref{MV36}) and (\ref{MV37}) respectively.

The proof of the theorem is trivial and employs Eqs. (\ref{MV31}). This theorem permits us to interpret the
functional
\begin{eqnarray}
S_{\kappa}=-\int_Adx \,\,y(x)\ln_{\kappa}y(x) \ ,
\label{MV40}
\end{eqnarray}
which can be written also in the form
\begin{eqnarray}
S_{\kappa}=\int_Adx \,\,\,\frac{y(x)^{1-\kappa}-y(x)^{1+\kappa}}{2\kappa} \ ,
\label{MV41}
\end{eqnarray}
as the entropy associated to the function
$\exp_{\kappa}\left(x\right)$. It is remarkable that in the $\kappa
\rightarrow 0$ limit, as $\ln_{\kappa}(y)$ and
$\exp_{\kappa}\left(x\right)$ approach $\ln (y)$ and $\exp(x)$
respectively, the new entropy reduces to the old Boltzmann-Shannon
entropy.

It is shown that the entropy $S_{\kappa}$ has the standard properties of Boltzmann-Shannon entropy:
is thermodynamically stable, is Lesche stable, obeys the Khinchin axioms of continuity,
maximality, expandability and generalized additivity.

\sect{$\kappa$-Trigonometry}

\subsection{$\kappa$-Hyperbolic Trigonometry}

The $\kappa$-hyperbolic trigonometry
can be introduced by defining the $\kappa$-hyperbolic sine and
$\kappa$-hyperbolic cosine
\begin{eqnarray}
\sinh_{\kappa}(x) =\frac{\exp_{\kappa}(x) -\exp_{\kappa}(-x)}{2} \
\ , \label{MVI1} \\
\cosh_{\kappa}(x) =\frac{\exp_{\kappa}(x) +\exp_{\kappa}(-x)}{2} \
\ , \label{MVI2}
\end{eqnarray}
starting from the $\kappa$-Euler formula
\begin{eqnarray}
\exp_{\kappa}(\pm x)=\cosh_{\kappa}(x)\pm \sinh_{\kappa}(x) \ \ .
\label{MVI3}
\end{eqnarray}

\begin{figure}[h] \centerline{
\includegraphics[width=.8\columnwidth,angle=0]
{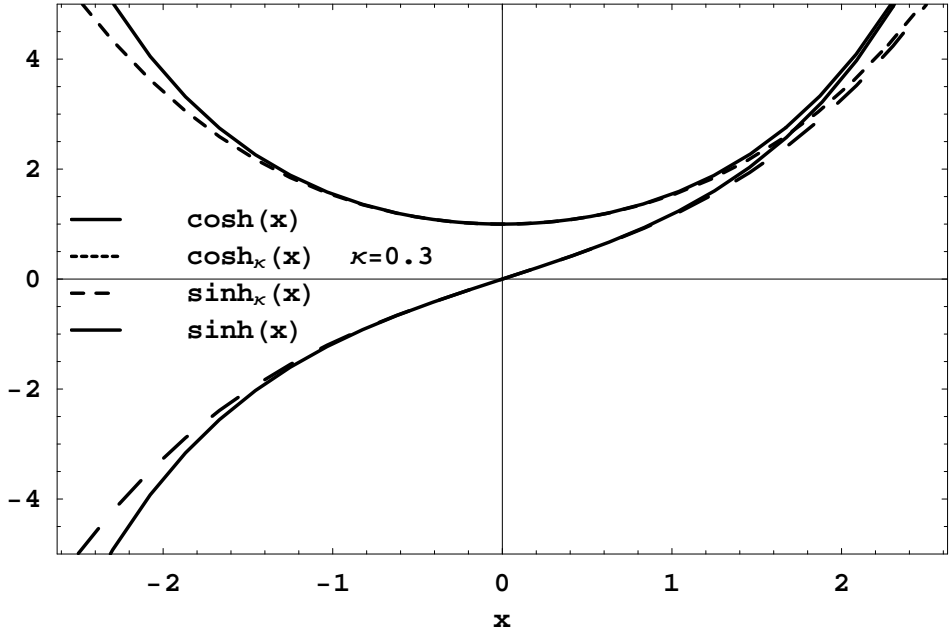}} \caption{Plot of  the functions $\sinh_{\kappa}(x)$ and
$\cosh_{\kappa}(x)$ for $\kappa=0.3$ (dashed lines) defined through
Eq. (\ref{MVI1}) and Eq. (\ref{MVI1}) respectively. For comparison in the same plot are reported the
corresponding ordinary functions $\sinh(x)$ and $\cosh(x)$ (continuous lines).} \label{FigA5}
\end{figure}

In Fig.\ref{FigA5} are plotted  the functions $\sinh_{\kappa}(x)$
and $\cosh_{\kappa}(x)$ for $\kappa=0.3$ (dashed lines). For
comparison, in the same figure are reported the corresponding
ordinary functions $\sinh(x)$ and $\cosh(x)$ (continuous lines).

The $\kappa$-hyperbolic tangent and cotangent functions are
defined through

\begin{eqnarray}
&&\tanh_{\kappa}(x)=\frac{\sinh_{\kappa}(x)}{\cosh_{\kappa}(x)} \
, \label{MVI4} \\
&&\coth_{\kappa}(x)=\frac{\cosh_{\kappa}(x)}{\sinh_{\kappa}(x)} \
. \label{MVI5}
\end{eqnarray}

Holding the relationships
\begin{eqnarray}
&&\sinh_{\kappa}(x)=\sinh\left(\{x\}\right) \ \ , \label{MVI6} \\
&&\cosh_{\kappa}(x)=\cosh\left(\{x\}\right) \ \ , \label{MVI7} \\
&&\tanh_{\kappa}(x)=\tanh\left(\{x\}\right) \ \ , \label{MVI8}  \\
&&\coth_{\kappa}(x)=\coth\left(\{x\}\right) \ \ , \label{MVI9}
\end{eqnarray}
it is straightforward to verify that $\kappa$-hyperbolic
trigonometry preserves the same structure of the ordinary
hyperbolic trigonometry which recovers as special case in the
limit $\kappa \rightarrow 0$. For instance from the $\kappa$-Euler formula and from $\exp_{\kappa}(-x)\exp_{\kappa}(x)=1$ the
fundamental formula of the $\kappa$-hyperbolic trigonometry
follows
\begin{eqnarray}
\cosh_{\kappa}^2(x)- \sinh_{\kappa}^2(x)=1 \ . \label{MVI10}
\end{eqnarray}

All the formulas of the ordinary hyperbolic trigonometry still hold,
after proper generalization. Taking into account that
$\{x\stackrel{\kappa}{\oplus}y \}=\{x\}+\{y\}$, it is easy to verify
that the generalization of a given formula can be obtained starting
from the corresponding ordinary formula, and then by making in the
arguments of the hyperbolic trigonometric functions the
substitutions $x+y\rightarrow x\stackrel{\kappa}{\oplus}y$ and
$x-y\rightarrow x\stackrel{\kappa}{\ominus}y$. For instance it
results
\begin{equation}
\sinh_{\kappa}(x\stackrel{\kappa}{\oplus}y )\!+\!\sinh_{\kappa}(
x\stackrel{\kappa}{\ominus}y ) = 2\sinh_{\kappa}(x)
\cosh_{\kappa}(y)\ \ , \ \ \ \label{MVI11}
\end{equation}
\begin{equation}
\cosh_{\kappa}(x\stackrel{\kappa}{\oplus}y
)=\cosh_{\kappa}(x)\cosh_{\kappa}(x)+\sinh_{\kappa}(x)
\sinh_{\kappa}(y)\ \ , \ \ \ \label{MVI12}
\end{equation}
\begin{eqnarray}
\tanh_{\kappa}(x)+\tanh_{\kappa}(y)
=\frac{\sinh_{\kappa}(x\stackrel{\kappa}{\oplus}y )}
{\cosh_{\kappa}(x)\cosh_{\kappa}(y)} \ \ , \ \ \ \label{MVI13}
\end{eqnarray}
and so on.

Obviously the substitution  $nx\rightarrow
x\stackrel{\kappa}{\oplus}x\stackrel{\kappa}{\oplus} ...
\stackrel{\kappa}{\oplus}x = [n]\stackrel{\kappa}{\otimes}
 x $ is required, so that, for instance it holds the formula
\begin{equation}
\sinh_{\kappa}^4(x)\!=\!\frac{1}{8}\!\left[\,
\cosh_{\kappa}\!\left([4]\stackrel{\kappa}{\otimes}
 x\right)\!-\!4\cosh_{\kappa}\!\left([2]\stackrel{\kappa}{\otimes}
 x\right)+3\right] , \label{MVI14}
\end{equation}
and so on.

The $\kappa$-De Moivre  formula involving hyperbolic
trigonometric functions having arguments of the type $r x$, with
$r\in{\bf R}$, assumes the form
\begin{equation}
[\cosh_{\kappa}(x)\pm \sinh_{\kappa}(x) ]^{r}
=\cosh_{{\kappa}/{\scriptstyle r}}(r
x)\!\pm\!\sinh_{_{\{{\scriptstyle \kappa}/{\scriptstyle
r}\}}}\!(r x) \ . \label{MVI15}
\end{equation}
Also the formulas involving the derivatives of the hyperbolic
trigonometric function still hold, after properly generalized. For
instance we have
\begin{eqnarray}
&&\frac{d \,\sinh_{\kappa}(x) }{d_{\kappa}x}= \cosh_{\kappa}(x) \ \ ,
\label{MVI16}  \\
&&\frac{d \,\tanh_{\kappa}(x) }{d_{\kappa}x}=
\,\cosh^{-2}_{\kappa}(x) \label{MVI17}  \ \ ,
\end{eqnarray}
and so on.

The $\kappa$-inverse hyperbolic functions can be
introduced starting from the corresponding ordinary functions. It is
trivial to verify that $\kappa$-inverse hyperbolic functions are
related to the $\kappa$-logarithm by the usual formulas of the
ordinary mathematics. For instance we have
\begin{eqnarray}
&&{\rm
arcsinh}_{\kappa}(x)=\ln_{\kappa}\left(\sqrt{1+x^2}+x\right)
\label{MVI18} \ \ , \\
&&{\rm
arccosh}_{\kappa}(x)=\ln_{\kappa}\left(\sqrt{x^2-1}+x\right)
\label{MVI19} \ \ , \\
&&{\rm arctanh}_{\kappa}(x)=\ln_{\kappa}\sqrt{\frac{1+x}{1-x}}
\label{MVI20} \ \ ,
\\
&&{\rm arccoth}_{\kappa}(x)=\ln_{\kappa}\sqrt{\frac{1-x}{1+x}}
\label{MVI21} \ \ ,
\end{eqnarray}
and consequently  hold
\begin{eqnarray}
{\rm arcsinh}_{\kappa}(x)&=&{\rm arccosh}_{\kappa}\sqrt{1+x^2} \ \ ,
\label{MVI22} \\
{\rm arcsinh}_{\kappa}(x)&=& {\rm arctanh}_{\kappa}\frac{x}{\sqrt{1+x^2}} \ \ ,
\label{MVI23} \\
{\rm arcsinh}_{\kappa}(x)&=&  {\rm
arccoth}_{\kappa}\frac{\sqrt{1+x^2}}{x}  \  \ . \label{MVI24}
\end{eqnarray}

From Eq. (\ref{MVI18}) follows  the relationship
\begin{eqnarray}
\exp_{\kappa}\left ({\rm arcsinh}_{\kappa}\, x\right)= \exp\left
({\rm arcsinh}\, x\right) \ \ . \label{MVI25}
\end{eqnarray}
Also the relationship
\begin{eqnarray}
{\rm arcsinh}_{\kappa}\left(x\right)=\frac{1}{\kappa}\,
\sinh_{\,1/\kappa}\left(\kappa x\right) \ \ , \label{MVI26}
\end{eqnarray}
involving the function ${\rm arcsinh}_{\kappa}(x)$, follows from Eq.
(\ref{MVI18}).  Analogous formulas involving ${\rm
arccosh}_{\kappa}(x)$, ${\rm arctanh}_{\kappa}(x)$ or ${\rm
arccoth}_{\kappa}(x)$ do not hold instead.

\subsection{$\kappa$-Cyclic Trigonometry}

By employing the generalized $\kappa$-Euler formula
\begin{eqnarray}
\exp_{\kappa}(\pm ix)=\cos_{\kappa}(x)\pm i\sin_{\kappa}(x) \ \ ,
\label{MVII1}
\end{eqnarray}
we introduce the $\kappa$-cyclic sine and $\kappa$-cosine as
\begin{eqnarray}
\sin_{\kappa}(x) =\frac{\exp_{\kappa}(ix) -\exp_{\kappa}(-ix)}{2i} \
\ , \label{MVII2} \\
\cos_{\kappa}(x) =\frac{\exp_{\kappa}(ix) +\exp_{\kappa}(-ix)}{2} \
\ , \label{MVII3}
\end{eqnarray}
while the $\kappa$-cyclic tangent and $\kappa$-cotangent functions are
defined through
\begin{eqnarray}
&&\tan_{\kappa}(x)=\frac{\sin_{\kappa}(x)}{\cos_{\kappa}(x)} \
, \label{MVII4} \\
&&\cot_{\kappa}(x)=\frac{\cos_{\kappa}(x)}{\sin_{\kappa}(x)} \ .
\label{MVII5}
\end{eqnarray}
After noting that
\begin{eqnarray}
\exp_{\kappa}(ix) = \exp\left(i\{x\}\right)  \ \ , \label{MVII6}
\end{eqnarray}
with
\begin{eqnarray}
\{x\}=\frac{1}{\kappa}\arcsin \kappa x  \ \ , \label{MVII7}
\end{eqnarray}
it follows that the cyclic functions are defined in the interval $-1/\kappa\leq x \leq 1/\kappa$. The function
\begin{eqnarray}
[x]=\frac{1}{\kappa}\sin \kappa x  \ \ , \label{MVII15}
\end{eqnarray}
is defined as the inverse of $\{x\}$, i.e $[\{x\}]=\{[x]\}=x$. The
$\kappa$-sum $\stackrel{\kappa}{\oplus}$ and $\kappa$-product
$\stackrel{\kappa}{\otimes}$ given by
\begin{eqnarray}
&&x\stackrel{\kappa}{\oplus}y=x\sqrt{1-\kappa^2y^2}+y\sqrt{1-\kappa^2x^2}
\ \  , \label{MVII16}
\\
&&x\stackrel{\kappa}{\otimes}y={1\over\kappa}\,\sin
\left(\,{1\over\kappa}\,\,{\rm arcsin}\,(\kappa x)\,\,{\rm
arcsin}\,(\kappa y)\,\right) \ , \label{MVII17}
\end{eqnarray}
are isomorphic operations to the ordinary sum and product respectively, i.e.
\begin{eqnarray}
&&\{x\stackrel{\kappa}{\oplus}y\}=\{x\}+\{y\} \ , \label{MVII18}
\\ &&\{x\stackrel{\kappa}{\otimes}y\}=\{x\}\cdot\{y\}  \ . \label{MVII19}
\end{eqnarray}

Holding the relationships
\begin{eqnarray}
&&\sin_{\kappa}(x)=\sin\left(\{x\}\right) \ \ , \label{MVII8} \\
&&\cos_{\kappa}(x)=\cos\left(\{x\}\right) \ \ , \label{MVII9} \\
&&\tan_{\kappa}(x)=\tan\left(\{x\}\right) \ \ , \label{MVII10}  \\
&&\cot_{\kappa}(x)=\cot\left(\{x\}\right) \ \ , \label{MVII11}
\end{eqnarray}
it is straightforward to verify that the generalized cyclic
trigonometry preserves the same structure of the ordinary cyclic
trigonometry, which recovers as special case in the limit $\kappa
\rightarrow 0$. For instance the following generalized formulas hold
\begin{eqnarray}
&&\cos_{\kappa}^2(x)+ \sin_{\kappa}^2(x)=1 \ , \label{MVII12}
\\
&&\sin_{\kappa}(x\stackrel{\kappa}{\oplus}y
)=\sin_{\kappa}(x)\cos_{\kappa}(y)+\cos_{\kappa}(x)
\sin_{\kappa}(y)\ \ , \ \ \ \label{MVII13}
\\
&&\cos_{\kappa}^5(x)\!=\!\frac{1}{16}\!\left[\,
\cos_{\kappa}\!\left([5]\stackrel{\kappa}{\otimes}
 x\right)\!+\!5\cos_{\kappa}\!\left([3]\stackrel{\kappa}{\otimes}
 x\right)+10\cos_{\kappa}(x)\right] , \label{MVII14}
\end{eqnarray}
and so on.

After introducing the following $\kappa$-deformed derivative operator
\begin{eqnarray}
\frac{d }{d _{\kappa}\,x}=\sqrt{1-\kappa^2\,x^2}\,\,\,\frac{d }{d \,
x} \ , \label{MVII20}
\end{eqnarray}
we can obtain easily further formulas oh the cyclic $\kappa$-trigonometry emerging as generalizations of the corresponding formulas of the ordinary trigonometry. For instance we have
\begin{eqnarray}
&&\frac{d \,\cos_{\kappa}(x) }{d_{\kappa}x}=- \sin_{\kappa}(x) \ \ ,
\label{MVII21}
\\
&&\frac{d \,\cot_{\kappa}(x) }{d_{\kappa}x}=-\,
\sin^{-2}_{\kappa}(x) \label{MVII22} \ \ ,
\end{eqnarray}
and so on.

The $\kappa$-inverse cyclic functions can be calculated by
inversion of the corresponding direct functions and are given by
\begin{eqnarray}
&&{\rm arcsin}_{\kappa}(x)=-i\,\ln_{\kappa}
\left(\sqrt{1-x^2}+ix\right) , \label{MVII23}
\\
&&{\rm
arccos}_{\kappa}(x)=-i\,\ln_{\kappa}\left(\sqrt{x^2-1}+x\right)
\label{MVII24} \ \ , \\
&&{\rm arctan}_{\kappa}(x)=i\,\ln_{\kappa}\sqrt{\frac{1-ix}{1+ix}}
\label{MVII25} \ \ ,
\\
&&{\rm arccot}_{\kappa}(x)=i\,\ln_{\kappa}\sqrt{\frac{ix+1}{ix-1}}
\label{MVII26} \ \ .
\end{eqnarray}

Finally we note that the $\kappa$-cyclic and $\kappa$-hyperbolic trigonometric functions are linked through the relationships
\begin{eqnarray}
&&\sin_{\kappa}(x)=-i\sinh_{\kappa}(ix) \ \ , \label{MVII27} \\
&&\cos_{\kappa}(x)=\cosh_{\kappa}(ix) \ \ , \label{MVII28} \\
&&\tan_{\kappa}(x)=-i\tanh_{\kappa}\left(ix\right) \ \ , \label{MVII29} \\
&&\cot_{\kappa}(x)=i\coth_{\kappa}\left(ix\right) \ \ , \label{MVII30} \\
&&{\rm arcsin}_{\kappa}(x)=-i\,{\rm arcsinh}_{\kappa}(ix)
\label{MVII31}  \ \ ,
\\
&&{\rm arccos}_{\kappa}(x)= -i\,{\rm arccosh}_{\kappa}(x)
\label{MVII32} \ \ , \\
&&{\rm arctan}_{\kappa}(x)=-i\,{\rm arctanh}_{\kappa}(ix)
\label{MVII33}  \ \ ,
\\
&&{\rm arccot}_{\kappa}(x)=i\,{\rm arccoth}_{\kappa}(ix)
\label{MVII34} \ \ ,
\end{eqnarray}
which in the $\kappa \rightarrow 0$ limit reduce to the standard
formulas involving the ordinary cyclic and hyperbolic functions.

\end{document}